\documentclass[12pt]{amsart}
\usepackage{amssymb,amsmath,amscd,enumerate}
\usepackage{tikz}
\usepackage{capt-of}
\usepackage[colorlinks=true,linkcolor=blue,citecolor=blue]{hyperref}
\usepackage{graphicx}
\usepackage{float} 
\usepackage{mathrsfs}
\usetikzlibrary{arrows}
 \numberwithin{equation}{section}
 \newtheorem{thm}{\bf Theorem}[section]
 \newtheorem{lem}[thm]{\bf Lemma}

 \newtheorem{prop}[thm]{\bf Proposition}
 \theoremstyle{definition}
 \newtheorem{defn}[thm]{Definition}

\newtheorem*{thm*}{Theorem}
\usepackage{color}

\DeclareMathOperator{\supp}{supp}
\DeclareMathOperator{\reg}{reg}
\DeclareMathOperator{\im}{im}
\DeclareMathOperator{\Inter}{Inter}

\def\a {\mathbf a}
\def\b {\mathbf b}
\def\c {\mathbf c}
\def\e {\mathbf e}
\def\m{\mathfrak m}
\def\NN{\mathbb N}
\newcommand {\kk} {\Bbbk} 

\usepackage[normalem]{ulem}

\begin{document}
\title[Powers of edge ideals of cubic circulant  graphs]{Regularity of powers and symbolic powers of edge ideals of cubic circulant graphs}
\date{}

\author[Nguyen Thu Hang]{Nguyen Thu Hang}
\address{Thai Nguyen University of Sciences, Tan Thinh Ward, Thai Nguyen City, Thai Nguyen, Vietnam}
\email{hangnt@tnus.edu.vn}

\author[My Hanh Pham]{My Hanh Pham}
\address{Faculty of Education, An Giang University\\Vietnam National University Ho Chi Minh City\\Dong Xuyen, Long Xuyen, An Giang, Vietnam}
\email{pmhanh@agu.edu.vn}

\author[Thanh Vu]{Thanh Vu}
\address{Institute of Mathematics, VAST, 18 Hoang Quoc Viet, Hanoi, Vietnam}
\email{vuqthanh@gmail.com}

\let\thefootnote\relax\footnotetext{*Corresponding Author: }
\keywords{Castelnuovo--Mumford regularity, cubic circulant graph, edge ideal, power and symbolic power, Minh's conjecture}

\subjclass[2000]{Primary: 05A05, 05A10; Secondary: 05C69}
 
\begin{abstract} 
 We compute the regularity of powers and symbolic powers of edge ideals of all cubic circulant graphs. In particular, we establish Minh's conjecture for cubic circulant graphs.
 \end{abstract}

\maketitle

\section{Introduction}
  \label{sec1} 
Let $G$ be a simple graph on the vertex set $V(G) = \{1,\ldots,n\}$ and edge set $E(G) \subseteq V(G) \times V(G)$. Let $R= \kk[x_1,\ldots, x_n]$ be a polynomial ring over a field $\kk$. The edge ideal of $G$ over $\kk$ is denoted by $I(G)={\big(} x_ix_j\mid\{i,j\}\in E(G) {\big)} \subseteq R$. Studying the homological invariants of powers of $I(G)$ and their relation to the combinatorial invariants of $G$ has been an active area of research. It shows many beautiful connections between algebra and combinatorics. In recent years, much attention has been devoted to Minh's conjecture stating that the regularity of powers of $I(G)$ and symbolic powers of $I(G)$ are equal. By the work \cite{ABS, BBH, F1, F2, F3, JS, MV}, the conjecture is known for unicyclic graphs, chordal graphs, Cameron-Walker graphs, and graphs of dimension one. In this paper, we establish Minh's conjecture for cubic circulant graphs by computing the regularity of all powers and symbolic powers of their edge ideals.

Circulant graphs have wide applications in coding, networks, and designs \cite{CJ, K, M}. Algebraic studies of edge ideals of circulant graphs have been carried out in \cite{BH, HMMP, R, UV, VVW}. In particular, Uribe-Paczka and Van Tuyl \cite{UV} computed the regularity of edge ideals of all cubic circulant graphs. Recall that for a positive integer $n \ge 3$ and a subset $S \subseteq \{1,\ldots,\lfloor n/2 \rfloor \}$, the circulant graph generated by $S$, denoted by $C_n(S)$, is the graph whose vertex set is $[n]$ in which two distinct vertices $i$ and $j$ are adjacent if and only if either $|i-j|\in S$ or $n-|i-j|\in S.$ Circulant graphs include complete graphs, cycle graphs, and crown graphs. Circulant graph $C_n(S)$ is regular of degree $2|S|$ unless $n$ is even and $n/2 \in S$, in this case it is regular of degree $2|S| - 1$. In particular, a cubic circulant graph is a graph of the form $C_{2n}(S)$ with $S = \{a,n\}$ for some integer $a$ such that $1 \le a < n$. For simplicity, we also denote it by $C_{2n}(a,n)$. Davis and Domke \cite{DD} proved that cubic circulant graphs are isomorphic to disjoint unions of certain copies of circulant graphs of the form $C_{2m}(1,m)$ and $C_{2m}(2,m)$ for $m > 2$ is an odd number. By the results of \cite{HNTT, NV} we may assume that $G$ is connected when studying the regularity of powers and symbolic powers of its edge ideal. For connected cubic circulant graphs, we prove:
 
\begin{thm}\label{thm_main} Let $G = C_{2n}(1,n)$ or $G= C_{2n}(2,n)$ where $n > 1$ is an odd number. Then 
$$\reg \left ( I(G)^t \right ) = \reg \left ( I(G)^{(t)} \right ) = 2t- 1 + \left \lfloor \frac{n}{2} \right \rfloor,$$
for all $t \ge 2$, where $I(G)^{(t)}$ denotes the $t$-th symbolic power of $I(G)$.
\end{thm}
 
In particular, the regularity stabilization index of edge ideals of connected cubic circulant graphs is either one or two depending on whether $\reg (I(G)) = \im(G) + 1$ or $\im(G) + 2$, where $\im(G)$ is the induced matching number of $G$. See Section \ref{sec2} for more details. By \cite[Theorem 4.5]{BHT} and \cite[Theorem 4.6]{GHOS}, it is known that 
\begin{equation}
    2t -1 + \im(G) \le \min \left \{ \reg \left ( I(G)^t  \right ) , \reg \left ( I(G)^{(t)} \right )  \right \}.
\end{equation}
By \cite[Theorem 3.3]{Ro}, we deduce that $\im(G) = \lfloor n/2\rfloor$, it remains to establish the upper bound. For regular powers, by induction and \cite[Theorem 5.2]{B}, it suffices to prove
\begin{equation}\label{eq1}
   \reg \left ( I(G)^t : (e_1 \cdots e_{t-1}) \right ) \le \im(G) + 1,
\end{equation}
for every tuple $\e = (e_1,\ldots,e_{t-1})$ of $t-1$ (possibly) repeated edges of $G$, where by abuse of notation, $e_j$ also denotes the minimal generator of $I(G)$ corresponding to the edge $e_j$ of $G$. By \cite[Theorem 6.7]{B}, we have that $I(G)^t : (e_1 \cdots e_{t-1}) = I(G_\e)$ where $G_\e$ is the graph obtained from $G$ by adding $\e$-even connection edges (the new edge might be a loop). Assume that $e_1 = \{u,v\}$ and $\e' = (e_2,\ldots,e_{t-1})$. To establish \eqref{eq1}, we prove the following 
\begin{enumerate}
    \item $\reg (I(G_\e)) \le \max \left \{ 2,  \reg \left ( \sqrt{I(G_\e)} \right ) \right \}$, 
    \item $\sqrt{I(G_\e) } = \sqrt{ I(G_{\e'}) : x_u} \bigcap \sqrt{ I(G_{\e'}) : x_v}$,
    \item $\reg \left ( \sqrt{ I(G_{\e'}) : x_u} \bigcap \sqrt{ I(G_{\e'}) : x_v} \right ) \le \reg (I(G_{\e'}))$.
\end{enumerate}
The conclusion then follows from induction on $t$. See Section \ref{sec3} for more details. For symbolic powers, we prove that for all exponents $x^\a = x_1^{a_1} \cdots x_n^{a_n} \in R$ such that $x^\a \notin I^{(t)}$ then 
$$\sqrt{I(G)^t : x^\a} = \sqrt{I(G)^{(t)} : x^\a}.$$ 
By \cite[Lemma 2.18]{MNPTV}, we deduce that $\reg \left ( I(G)^{(t)}  \right ) \le \reg (I(G)^t).$

In the next section, we recall some properties of cubic circulant graphs and their edge ideals. In Section \ref{sec3}, we compute the regularity of powers of edge ideals of connected cubic circulant graphs. In Section \ref{sec4}, we compute the regularity of symbolic powers of edge ideals of connected cubic circulant graphs. We then deduce the formula for the regularity of powers and symbolic powers of edge ideals of arbitrary cubic circulant graphs.

\section{Regularity of edge ideals of cubic circulant graphs}\label{sec2}
In this section, we recall relevant concepts concerning graphs, in particular circulant graphs and the regularity of their edge ideals.

 \subsection{Regularity} Let $R = \kk[x_1,\ldots,x_n]$ be a polynomial ring over a field $\kk$. Denote by $\m = (x_1, \ldots, x_n)$ the maximal homogeneous ideal of $R$. For a finitely generated graded $R$-module $M$, the Castelnuovo-Mumford regularity (or regularity for short) of $M$ is defined by 
$$\reg_R(M) = \max \left \{ i + j \mid H^i_\m (M)_j \neq 0 \right \},$$
where $H^i_\m (M)$ denotes the $i$th local cohomology of $M$ with respect to $\m$. For a non-zero homogeneous ideal $I$ of $R$, we have 
      $$\reg_R(I) = \reg_R(R/I) +1.$$

 \subsection{Graph} 
 
Let $G$ denote a finite simple graph over the vertex set $V(G) = [n] = \{1,\ldots,n\}$ and the edge set $E(G)$. A simple graph $H$ is a subgraph of $G$ if $V(H) \subseteq V(G)$ and $E(H) \subseteq E(G)$. $H$ is an induced subgraph of $G$ if it is a subgraph of $G$ and $E(H) = E(G) \cap V(H) \times V(H)$. $G$ is called bipartite if there exists a partition $V(G) = X \cup Y$ such that $E(G) \subseteq X \times Y$.

A cycle $C_n$ of length $n \ge 3$ is a graph on $[n]$ whose edges are $\{i,i+1\}$ for $i = 1, \ldots, n-1$ and $\{1,n\}$. The complement of $G$, denoted by $\overline{G}$, is the graph whose vertex set is $V(G)$ and the edge set consists of all edges $\{i,j\}$ such that $\{i,j\}\notin E(G).$ $G$ is called weakly chordal if every induced cycle in both $G$ and its complement $\overline{G}$ has length at most $4$.

For any graph $G$, a subset $M \subseteq E(G)$ is called a \textit{matching} of $G$ if no two
edges in $M$ share a common vertex. It is an \textit{induced matching} if $M$ forms an induced subgraph of $G$. The size of a maximum induced matching of $G$, denoted by $\im(G)$, is called the\textit{ induced matching number} of $G$.

For a vertex $i\in V(G)$, let the neighborhood of $i$ be the subset $N_G(i)=\{j\in V(G) \mid \{i,j\}\in E(G)\}$. The closed neighborhood of $i$ is $N_G[i]=N_G(i)\cup\{i\}$. Let $U$ be a subset of the vertex set of $G$. We denote by $G[U]$, respectively $G\backslash U$ the induced subgraph of $G$ on $U$, respectively on $V(G) \backslash U$. When $U=\{i\}$, we write $G\backslash i$ instead of $G\backslash \{i\}.$
 The neighborhood of $U$ in $G$ is the set
\[
N_{G}(U):=\{v\in V(G) \mid \{u,v\}\in E(G) \text{ for some } u\in U\}.
\]
The closed neighborhood of $U$ is $N_{G}[U]:=U\cup N_{G}(U)$.

\medskip

Let $1 \le a < n$ be positive integers. The cubic circulant graph $G = C_{2n}(a,n)$ is a graph on vertex set $V(G)  = [2n]$ and $\{i,j\}$ is an edge of $G$ if and only if $|i-j| \in \{a,n,2n-a\}$. It is easy to see that $G$ is connected if and only if $\gcd(a,n) = 1$. Davis and Domke \cite{DD} proved the following structural results for cubic circulant graphs.
\begin{lem}\label{lem_structure_cubic}
	Let $G=C_{2n}(a,n)$ with $1\leq a <n$ and let $t=\gcd(a,2n).$
	\begin{itemize}
		\item [(i)] If $\frac{2n}{t}$ is even, then $G$ is isomorphic to $t$ copies of  $C_\frac{2n}{t}(1,\frac{n}{t}).$
		\item [(ii)] If $\frac{2n}{t}$ is odd, then $G$ is isomorphic to $\frac{t}{2}$ copies of $C_{\frac{4n}{t}}(2, \frac{2n}{t}).$
	\end{itemize}
\end{lem}

From \cite[Theorem 3.3]{Ro}, we deduce an explicit formula for the induced matching numbers of connected cubic circulant graphs.

 \begin{lem}\label{lem_induced_matching_cubic}
 	Let $G = C_{2n}(a,n)$ for $a \in \{1,2\}$ be a connected cubic circulant graph. Then
  $$ \im(G)= \left \lfloor
        \frac{n}{2} \right \rfloor.$$
\end{lem}
 \begin{proof}
 With the terminology as in \cite[Theorem 3.3]{Ro}, we have $|E(G)| = 3n$, $s = 2$, and $|A| = 1$. Hence, $$\im(G) = \left \lfloor \frac{|E(G)|}{s^2 + (|A| + 1) s - 2} \right \rfloor = \left \lfloor \frac{n}{2}  \right \rfloor.$$
 The conclusion follows.
  \end{proof}

Finally, we have the following result on the odd cycles of cubic circulant graphs that will be useful later. 

\begin{lem}\label{lem_cycle_cubic} Let $G = C_{2n}(a,n)$ for $a \in \{1,2\}$ be a connected cubic circulant graph. Let $\Delta$ be an odd cycle of $G$. Then $N_G(\Delta) = V(G)$.
\end{lem}
 \begin{proof} Note that $G \setminus \{i,i+n\}$ is bipartite for any $i = 1, \ldots, n$. Hence, $V(\Delta) \cap \{i,i+n\} \neq \emptyset$ for all $i = 1, \ldots, n$. The conclusion follows.   
 \end{proof}

\subsection{Edge ideals of graphs and their regularity} Let $G$ be a simple graph on the vertex set $V(G) = [n]$ and the edge set $E(G)$. The edge ideal of $G$ is defined by 
$$I(G) = (x_ix_j \mid \{i,j\} \in E(G)) \subseteq R.$$
We have the following result on the regularity of edge ideals of connected cubic circulant graphs. 
\begin{prop}\label{prop1} Let $G = C_{2n}(a,n)$ for $a \in \{1,2\}$ be a connected cubic circulant graph. If $n \le 4$ then $\reg (I(G)) = \im(G) + 1$. If $n \ge 5$, then 
	$$\reg(I(G)) = \begin{cases}
\im(G) + 2, & \text{ if } a = 1 \text{ and } n= 4k + 1 \text{  or if } a =2 \text{ and } n = 4k+3, \\
     \im(G)+1, & \text{ otherwise}.
    \end{cases}$$
\end{prop}
\begin{proof}
    The conclusion follows from Lemma \ref{lem_induced_matching_cubic} and \cite[Lemma 4.3]{UV}.
\end{proof}


\section{Regularity of powers of edge ideals of cubic circulant graphs}\label{sec3}
In this section, we compute the regularity of powers of edge ideals of connected cubic circulant graphs. We will use the result of Banerjee \cite{B} to bound the regularity of powers of edge ideals in terms of the regularity of colon ideals by a product of edges of $G$. We first analyze these colon ideals in more detail.
\subsection{Colon ideals of powers of edge ideals by products of edges} In this subsection, we let $G$ be an arbitrary simple graph. Let $\e = (e_1, \ldots, e_{t-1})$ be a tuple of $t-1$ possibly repeated edges of $G$. 

\begin{defn}
	Two vertices $u$ and $v$ ($u$ may be the same as $v$) are said to be \textit{even-connected} with respect to a tuple $\e = (e_1, \ldots, e_{t-1})$ if there is a path $p_0p_1\dots p_{2k+1}$, $k\ge 1$ in $G$ such that:
	\begin{enumerate}
		\item $p_0=u$ and $p_{2k+1}=v$.
		\item For all $0\le l\le k-1,$ $p_{2l+1}p_{2l+2} = e_i$ for some $i$.
		\item For all $i$, $$|\{l\ge 0\mid  p_{2l+1}p_{2l+2}=e_i\}|\le |\{j\mid e_j=e_i\}|.$$
		\item For all $0\le r \le 2k,$ $\{p_r,p_{r+1}\}$ is an edge in $G$.
	\end{enumerate}
\end{defn}
We denote by $G_\e$ the graph on the vertex set $V(G)$ and edge set 
$$E(G_\e) = E(G) \cup \{ \{u, v\}  \mid u,v \text{ are even-connected with respect to } \e\}.$$
Note that $G_\e$ might contain a loop. By abuse of notation, for an edge $e$ of $G$, we also use $e$ to denote the minimal generator of $I(G)$. It should be clear from the context when $e$ corresponds to an edge or a monomial in $R = \kk[x_1,\ldots,x_n]$. We denote by $x^\e = e_1 \cdots e_{t-1}$ the product of edges of $G$. With this notation, \cite[Theorem 6.7]{B} can be written as 
\begin{thm}\label{thm_colon_Ban} Let $G$ be a simple graph and $\e = (e_1, \ldots, e_{t-1})$ be a tuple of $t-1$ edges of $G$. Then 
$$I(G)^t : x^\e = I(G_\e).$$    
\end{thm}
The following lemma is a generalization of \cite[Lemma 3.4]{AB}.
\begin{lem}\label{lem_squarefree_colon} Let $G$ be a simple graph and $\e = (e_1^{a_1},\ldots,e_s^{a_s})$ be a tuple of edges of $G$ where $e_1, \ldots, e_s$ are distinct edges of $G$ and $a_1, \ldots, a_s > 0$ are the multiplicities of these edges in the tuple $\e$. Denote by $t-1 = a_1 + \cdots + a_s$. Assume that $I(G)^t : x^\e$ is squarefree. Then 
$$I(G)^t : x^\e = I(G)^{s+1} : (e_1 \cdots e_s) = \left ( I(G)^2 : e_1 \right )^{s} : (e_2 \cdots e_{s}).$$    
\end{lem}
\begin{proof}We denote by $I = I(G)$. We first prove that 
\begin{equation}\label{eq_sq_cl_1}
    I^t : x^\e = I^{s+1} : (e_1 \cdots e_s).
\end{equation}
Clearly, the left-hand side contains the right-hand side. Now assume that $u,v$ are $\e$-even connected. We prove Eq. \eqref{eq_sq_cl_1} by induction on $|\a|$. The base case where $a_1 = \cdots = a_s = 1$ is vacuous. Let $p_0\ldots p_{2k+1}$ be an even walk connecting $u$ and $v$. We may assume that there exist $0\le a < b < k$ such that $e_1 = \{p_{2a+1}, p_{2a+2} \} = \{p_{2b+1},p_{2b+2} \}$. In particular, the multiplicity $a_1$ of $e_1$ is greater than $1$. If $p_{2a+1} = p_{2b+1}$, then we have a shorter even walk $p_0\ldots p_{2a} p_{2b+1} p_{2b+2} \ldots p_{2k+1}$ connecting $u$ and $v$. In particular, $uv \in I^{t-1}: x^{\e'}$ where $\e' = (e_1^{a_1-1},\ldots,e_s^{a_s})$. By induction, we deduce that $uv \in I^{s+1} : (e_1 \cdots e_s)$. If $p_{2a+2} = p_{2b+1}$ then we have an even-closed walk connecting this repeated vertex. By Theorem \ref{thm_colon_Ban}, we deduce that $I^t:x^\e$ is not squarefree, a contradiction. 

We now prove
\begin{equation}\label{eq_sq_cl_2}
    I^{s+1} : (e_1 \cdots e_s) \subseteq (I^2 : e_1)^{s} : (e_2 \cdots e_{s}).
\end{equation}
Let $uv$ be a minimal generator of $I^{s+1} : (e_1 \cdots e_s)$ and $p_0 \ldots p_{2k+1}$ be an even walk connecting $u$ and $v$. If $e_1$ does not appear in the edges $\{p_{2l+1}, p_{2l+2}\}$ for any $l$ then $uv \in I^s : (e_2 \cdots e_s) \subseteq (I^2 : e_1)^s : (e_2 \cdots e_s)$. Hence, we may assume that $e_1 = \{p_{2l+1}, p_{2l+2}\}$. In particular, $p_{2l} p_{2l+3} \in I^2 : e_1$. By Theorem \ref{thm_colon_Ban}, we deduce that $uv \in (I^2 : e_1)^s : (e_2 \cdots e_s)$. 

It remains to prove 
\begin{equation}\label{eq_sq_cl_3}
    I^{s+1} : (e_1 \cdots e_s) \supseteq (I^2 : e_1)^{s} : (e_2 \cdots e_{s}).
\end{equation}
Let $I(G') = I^2 : e_1$ and $uv \in I(G')^s : (e_2 \cdots e_s)$. By Theorem \ref{thm_colon_Ban}, there exists an even walk in $G'$ with respect to $e_2,\ldots, e_s$ connecting $u$ and $v$. Now, for any edge in $G'$ which is not in $G$, we can replace with an even walk in $G$ via $e_1$. Concatenating these walks, we have an even walk in $G$ with respect to $e_1^{a_1}, e_2, \ldots, e_s$ where $a_1$ is the number of times we have an edge in $G'$ which is not in $G$. By Theorem \ref{thm_colon_Ban} and Eq. \eqref{eq_sq_cl_1}, we deduce that 
$$uv \in I^{a_1 + s} : (e_1^{a_1} e_2 \cdots e_s) \subseteq I^{s+1}:(e_1 \cdots e_s).$$
The conclusion follows.
\end{proof}

\subsection{Colon ideals of powers of edge ideals of cubic circulant graphs} In this subsection, we analyze in more details the colon ideals of powers of edge ideals of connected cubic circulant graphs.

\begin{lem}\label{lem_power_2} Let $G = C_{2n}(a,n)$ for $a \in \{1,2\}$ be a connected cubic circulant graph. Let $e$ be any edge of $G$. Then 
$$ \reg (I(G)^2 : x^e) \le \im(G) + 1.$$ 
\end{lem}
\begin{proof} If $n \le 3$ the conclusion follows by direct inspection. Thus, we assume that $n \ge 4$ and $e = \{ u,v\}$. By Lemma \ref{lem_cycle_cubic} and Theorem \ref{thm_colon_Ban} we deduce that $I(G)^2 : x^e$ is squarefree. In particular, 
 \begin{equation}\label{eq_power_2}
     I(G_e) = I(G)^2 : x^e = (I(G) : x_u) \cap (I(G) : x_v).
 \end{equation}
 Now, note that 
 \begin{align*}
     I(G) : x_u & = I(G \backslash N_G[u]) + (x_j \mid j \in N_G(u)) \\
     I(G) : x_u + I(G) : x_v & = I(G \backslash N_G[\{u,v\}] ) + ( x_j \mid j \in N_G[\{u,v\}]).
 \end{align*}
 Let $M$ be any induced matching of $G \backslash N_G[\{u,v\}]$ then $M \cup \{e\}$ is an induced matching of $G$. Hence, $\im(G \backslash N_G[\{u,v\}]) + 1 \le \im(G)$. Furthermore, for any $u \in V(G)$, $G \backslash N_G[u]$ is weakly chordal. By Eq. \eqref{eq_power_2} and \cite[Theorem 14]{W}, we deduce that 
 \begin{align*}
     \reg (I(G_\e)) &\le \max \left \{\reg  \left ( I(G \backslash N_G[u]) \right ), \reg \left ( I(G\backslash N_G[v]) \right ) , \reg  \left ( I(G\backslash N_G[\{u,v\}]) \right ) + 1 \right \} \\
     & = \max \left \{ \im \left ( G \backslash N_G[u] \right ) + 1, \im \left (G\backslash N_G[v] \right ) + 1 , \im (G\backslash N_G [ \{u,v\} ] ) + 2  \right \} \\
     & \le \im(G) + 1.
 \end{align*}
 The conclusion follows.
\end{proof}

We now have some preparation lemmas for the case that $G_\e$ contains a loop. 
\begin{lem}\label{lem_squarefree_multiplicity} Let $G = C_{2n}(a,n)$ for $a \in \{1,2\}$ be a connected cubic circulant graph. Let $\e = (e_1^{a_1},\ldots,e_s^{a_s})$ be a tuple of edges of $G$ with $|\a| = t-1$. Then 
$$ I(G)^t : x^\e = I(G)^{s+1} : (e_1 \cdots e_s).$$    
\end{lem}
\begin{proof}
    Let $u,v$ be $\e$-even connected. Assume that in the even walk $p_0 \ldots p_{2k+1}$ there are repeated edges, i.e. $\{p_{2a+1}, p_{2a+2}\} = \{p_{2b+1}, p_{2b+2} \}$ for some $a < b$. As in the proof of Lemma \ref{lem_squarefree_colon} if $p_{2a+1} = p_{2b+1}$ then we can shorten the walk and we are done by induction. Thus, we may assume that $p_{2a+1} = p_{2b+2}$. The closed even walk from $p_{2a+2}$ to $p_{2b+1}$ decomposes into cycles of $G$, at least one of them is odd, called $\Delta$. Now, by Lemma \ref{lem_cycle_cubic}, we deduce that $p_{2k+1} \in N_G(\Delta)$. Assume that $p_{2k+1} \in N_G(p_c)$ for some $c$ such that $ 2a+2 \le c < 2b + 1$. If $c$ is even, then we walk to $p_c$ and then go to $p_{2k+1}$. If $c$ is odd, from $p_{2a+1}$ we walk back from $p_{2b+2}$ to $p_c$ then we go to $p_{2k+1}$. Hence, we have a proper even connection that connects $u$ and $v$. By induction, the conclusion follows.
\end{proof}

\begin{lem}\label{lem_colon_loops} Let $G = C_{2n}(a,n)$ for $a \in \{1,2\}$ be a connected cubic circulant graph. Let $\e = (e_1, \ldots, e_{t-1})$ be a tuple of edges of $G$. Assume that $G_\e$ has a loop at $j$. Then $\{j,v\} \in E(G_\e)$ for all $v \in V(G)$. 
\end{lem}

\begin{proof} We may assume that $G_\e$ has a loop at $1$. Let $p_0\ldots p_{2k+1}$ be a closed even walk with $p_0 = p_{2k+1} = 1$. The closed even walk $p_0\ldots p_{2k+1}$ must decompose into simple cycles in $G$, at least one of which is odd, say $\Delta$. By Lemma \ref{lem_cycle_cubic}, $N_G(\Delta) = V(G)$. Let $v \in V(G)$ be an arbitrary vertex. Then $v \in N(p_c)$ for some $c$. We prove that $\{p_0,v\} \in G_\e$. Indeed, if $c$ is even, then the path $p_0\ldots,p_c,v$ is an $\e$-even walk connecting $p_0$ and $v$. If $c$ is odd, the path $p_{2k+1}, \ldots,p_c,v$ is an $\e$-even walk connecting $p_{2k+1} = p_0$ to $v$. The conclusion follows.
\end{proof}

\begin{lem}\label{lem_colon_loops_2} Let $G = C_{2n}(a,n)$ for $a \in \{1,2\}$ be a connected cubic circulant graph. Let $\e = (e_1, \ldots, e_s)$ be a tuple of distinct edges of $G$. Assume that $e_{s} = \{u,v\}$. We denote by $\e' = (e_1, \ldots, e_{s-1})$. Then 
$$ \sqrt { I(G_\e) }  =  \sqrt { I(G_{e'}) : x_u } \bigcap \sqrt { I(G_{e'}) : x_v }.$$ 
\end{lem}
\begin{proof} We first claim that 
\begin{equation}\label{eq_radical_1}
I(G_\e)   \subseteq I(G_{e'}) : x_u.    
\end{equation}
Let $\{w,y\}$ be an edge of $G_\e$. If $\{w,y\} \in E(G)$, there is nothing to prove. By Theorem \ref{thm_colon_Ban}, we may assume that there is an even walk $p_0 \cdots p_{2k+1}$ connecting $w$ and $y$. If none of the edges $\{p_{2l+1}, p_{2l+2}\}$ is $e_s$ then $\{u,v \} \in E(G_{\e'})$ and we are done. Hence, we assume that $p_{2l+1} = u$. In other words, $\{w,u\} \in E(G_{\e'})$. In particular, $x_w \in I(G_{\e'}) : x_u$. Eq. \eqref{eq_radical_1} follows.

We now claim that 
\begin{equation}\label{eq_radical_2}
\sqrt { I(G_{e'}) : x_u } \bigcap \sqrt { I(G_{e'}) : x_v } \subseteq \sqrt { I(G_\e) }.
\end{equation}
Since $G_{\e'}$ is a subgraph of $G_\e$, it suffices to prove that if $x_w \in I(G_{e'}) : x_u$ and $x_y \in I(G_{e'}) : x_v$ then $x_w x_y \in I(G_\e)$. This follows from Lemma \ref{lem_squarefree_multiplicity}.    
\end{proof}

We also have two simple lemmas.
\begin{lem}\label{lem_colon_m} Let $I$ be a non-zero monomial ideal and $x$ is a variable of $R$ such that $I:x = \m$. Then $\reg (I) \le \max \left \{2, \reg (I + (x)) \right \}.$    
\end{lem}
\begin{proof}
    The conclusion follows from the standard regularity lemma applied to the short exact sequence 
    $$0 \to x(I:x) \to I \oplus (x) \to I + (x) \to 0,$$
    and the fact that $\reg (\m) = \reg (x) = 1$.
\end{proof}

\begin{lem}\label{lem_reg_inter} Assume that $G$ is a graph with possibly loops. Let $\{u,v\}$ be a non-loop edge of $G$. Then 
$$\reg \left ( \sqrt { I(G) : x_u} \cap \sqrt{ I(G) : x_v} \right ) \le \reg (I (G)).$$    
\end{lem}
\begin{proof}
    Let $J = \sqrt { I(G) : x_u}$, $K = \sqrt { I(G) : x_v}$, and $L = J \cap K$. Then we have 
    $$\reg (L) \le \max \{ \reg (J), \reg (K), \reg (J + K) + 1\}.$$
    Note that $\reg (J + K)  = \reg \left (  \sqrt{ I( G \setminus N[\{u,v\}] ) } \right ) = \reg (I(H) )$, where $H = G \backslash U$ where $U = N[\{u,v\}] \cup \{ j \mid j \text{ is a loop of } G\}$. Now, we have $J, K$, and $I (H \cup \{u,v\})$ are restrictions of $I(G)$. The conclusion follows from \cite[Corollary 2.22]{MNPTV}.
\end{proof}

We are now ready for the main result of this section.
\begin{thm}\label{thm_reg_power}
	Let $G = C_{2n}(a,n)$ for $a \in \{1,2\}$ be a connected cubic circulant graph. Then 
 $$\reg (I(G)^t) = 2t -1 + \im(G),$$ 
 for all $t\ge 2.$
\end{thm} 
\begin{proof}
By \cite[Theorem 4.5]{BHT} and \cite[Theorem 5.7]{B}, it suffices to prove the following inequality for any tuple $\e = (e_1, \ldots, e_{t-1})$ of edges of $G$ by induction on $t \ge 2$ 
		\begin{equation}\label{eq_reg_p_1}
		    \reg (I(G_\e)) = \reg (I^{t}:e_1\ldots e_{t-1})\leq \im(G)+1.
		\end{equation}
		The base case $t = 2$ is Lemma \ref{lem_power_2}. Now, assume that $t > 2$. There are two cases. 

  \smallskip
\noindent \textbf{Case 1.} $G_\e$ has no loop. The conclusion follows from induction, Lemma \ref{lem_squarefree_colon} and \cite[Theorem 3.1]{BN}.
\smallskip

\noindent \textbf{Case 2.} $G_\e$ has a loop. Assume that $j_1, \ldots, j_s$ are loops of $G_\e$. By Lemma \ref{lem_colon_loops} and Lemma \ref{lem_colon_m}, we deduce that 
$$\reg (I(G_\e) ) \le \max \left \{2, \reg I(G_\e \backslash \{j_1,\ldots,j_s\}) \right \}.$$
By Lemma \ref{lem_squarefree_multiplicity}, we may assume that $e_1, \ldots, e_{t-1}$ are distinct. Let $\e' = (e_1,\ldots,e_{t-2})$ and $e_{t-1} = \{u,v\}$. By Lemma \ref{lem_colon_loops_2} and Lemma \ref{lem_reg_inter}, we have
\begin{align*}
\reg \left (  I( G_\e \setminus \{j_1,\ldots, j_s\}) \right ) &= \reg \left  ( \sqrt {I(G_e)} \right ) \\
& = \reg \left ( \sqrt { I(G_{\e'}) : x_u} \cap \sqrt{ I(G_{\e'}) : x_v} \right )  \le \reg (I (G_{\e'}) ).
\end{align*}
That concludes the induction step and the proof of the theorem.
\end{proof}

\section{Regularity of symbolic powers of edge ideals of cubic circulant graphs}\label{sec4}
In this section, we prove the strong form of the Minh's conjecture \cite[Conjecture A]{MV} for connected cubic circulant graphs. We then deduce the formula for the regularity of powers and symbolic powers of edge ideals of general cubic circulant graphs. We first recall some definitions. 

Let $I \subseteq R = \kk[x_1,\ldots,x_n]$ be a squarefree monomial ideal with the primary decomposition $I = P_1 \cap \cdots \cap P_r$. The $s$-th symbolic power of $I$ is defined by $$I^{(s)} = P_1^s \cap \cdots \cap P_r^s.$$

\begin{defn} Let $J \subseteq K$ be monomial ideals in $R$. We denote $\Inter(J , K)$ the set
of intermediate ideals between $J$ and $K$ containing all monomial ideals $L$ such that $L = J + ( f_1,\ldots, f_t)$ where $f_j$ are among minimal monomial generators of $K$.   
\end{defn}

For an exponent $\a = (a_1,\ldots,a_n) \in \NN^n$, $x^\a$ denotes the monomial $x^\a = x_1^{a_1} \cdots x^{a_n}$. The support of $\a$ is denoted by $\supp \a = \{ j\mid a_j \neq 0\}$. For a subset $F \subseteq [n]$, we denote by $x_F = \prod_{j \in F} x_j$. For a monomial $f$, we denote by $\partial (f) / \partial(x^\a) = f / \gcd(f,x^\a)$. 
The following is the key lemma in computing the regularity of symbolic powers of edge ideals of cubic circulant graphs.

\begin{lem}\label{lem_radical_colon} Let $G$ be a connected cubic circulant graph and $\a$ be an exponent such that $x^\a \notin I(G)^{(t)}$. Then $\sqrt{I(G)^t: x^\a} = \sqrt{I(G)^{(t)} : x^\a}$. 
\end{lem}
\begin{proof} By \cite[Theorem 5.9]{SVV}, we may assume that $G$ is not bipartite. It suffices to prove that if $x_F \in \sqrt{I(G)^{(t)} : x^\a}$ then $x_F \in \sqrt{I(G)^t : x^\a}$. We may assume that $F$ is an independent set. Let $U = N_G(F)$ and $x^\b$ be the restriction of $x^\a$ on $U$. Since $F$ is an independent set $U \cap F = \emptyset$. If $|\b| \ge t$, then $x_F^t x^\b \in I^t$. Hence, we may assume that $|\b| < t$. Let $s$ be an exponent such that $x_F^s \in I(G)^{(t)} : x^\a$. By \cite[Lemma 2.6]{MNPTV}, we deduce that 
$$ \partial (x_F^s x^\a )/ \partial(x^\b) = x_F^s x^\c \in I(G)^{(t - |\b|)}$$ 
where $\c = \a - \b$. In particular, $\supp \c \cap U = \emptyset$. Since $F$ is an independent set, this implies that $x_F^s x^\c \in I(H)^{(t - |\b|)}$, where $H$ is the induced subgraph of $G$ on $V(G) \setminus N_G[F]$. Since $H$ is bipartite, by \cite[Theorem 5.9]{SVV}, we deduce that $\partial (x_F^s x^\a) / \partial (x^\b) \in I(H)^{t - |b|}$. The conclusion follows from \cite[Lemma 2.18]{MV}.
\end{proof}
We also have the following simple lemma.
\begin{lem}\label{lem_lower_bound}
    Let $G$ be a simple graph and $L \in \Inter \left ( I(G)^s, I(G)^{(s)} \right )$ for some $s \ge 1$. Then 
    $\reg (L) \ge 2t - 1 + \im(G)$.
\end{lem}
\begin{proof}
    Let $H$ be an induced subgraph of $G$ consisting of in $\im(G)$ matching of $G$. By \cite[Corollary 2.7]{MNPTV}, we deduce that $L_{V(H)} = I(H)^s$. The conclusion then follows from \cite[Corollary 2.22]{MNPTV}.
\end{proof}
 \begin{thm}\label{thm_intermediate} Let $G$ be a connected cubic circulant graph and $L$ be an intermediate ideal lying between $I(G)^t$ and $I(G)^{(t)}$ for some $t \ge 2$. Then $\reg (L) = 2t - 1 + \im(G)$.     
 \end{thm}
 \begin{proof}
     By Lemma \ref{lem_radical_colon}, we deduce that for any exponent $\a$ such that $x^\a \notin L$ then $\sqrt{(I(G)^t : x^\a} = \sqrt{L:x^\a}$. By \cite[Lemma 2.18]{MNPTV}, we have $\reg (L) \le \reg (I(G)^t)$. The conclusion then follows from Theorem \ref{thm_reg_power} and Lemma \ref{lem_lower_bound}.
 \end{proof}
\begin{proof}[Proof of Theorem \ref{thm_main}] The conclusion follows from Lemma \ref{lem_induced_matching_cubic} and Theorem \ref{thm_intermediate}.
\end{proof}

We will now derive a formula for the regularity of powers and symbolic powers of general cubic circulant graphs. We first need two simple lemmas.

\begin{lem}\label{lem_mixed_sum_1} Let $G = kG_1$ be the disjoint union of $k$ copies of $G_1$. Assume that $\reg (I(G_1)^s) = 2s - 2 + \reg (I(G_1))$ for all $s \ge 2$. Then 
$$\reg (I(G)^t) = \reg (I(G)^{(t)}) = \reg (I(G)) + 2t- 2 = k \reg (I(G_1)) - k + 2t - 1$$ 
for all $t \ge 2$.
\end{lem}
\begin{proof}
	The conclusion follows from induction on $k$, \cite[Theorem 1.1]{NV}, and \cite[Theorem 4.7]{HNTT}.
\end{proof}

\begin{lem}\label{lem_mixed_sum_2} Let $G = kG_1$ be the disjoint union of $k$ copies of $G_1$. Assume that $\reg (I(G_1)^s) = 2s - 4 + \reg (I(G_1)^2)$ for all $s \ge 3$ and $\reg (I(G_1)^2) = \reg (I(G_1)) + 1$. Then 
$$\reg (I(G)^t)  = \reg (I(G)^{(t)}) = \reg (I(G)) + 2t -3 = k \reg (I(G_1)) -k + 2t -2$$
for all $t \ge 2$.
\end{lem}
\begin{proof}
	For completeness, we prove the equality for regular powers. The base case $k = 1$ follows from the assumption. Let $H = (k-1) G_1$ be the disjoint union of $k-1$ copies of $G_1$. Let $P = I(G_1)$, $Q = I(H)$, and $I = I(G) = P + Q$. By induction, we have 
 $$  \reg (Q^{j+1}) - \reg (Q^j) =   \reg (P^{j+1}) - \reg (P^j) = \begin{cases} 1 & \text { if } j = 1, \\
 2 & \text { if } j \ge 2.\end{cases}$$
 In other words, 
 \begin{equation}\label{eq_mixed_sum}
     \max_{i \in [1,t-1]} \left \{ \reg ( P ^{t-i}) + \reg (Q^i) \right \} = \reg (P^{t-1} ) + \reg ( Q).
 \end{equation}
Assume that $t \ge 2$. By \cite[Theorem 1.1]{NV} and Eq. \eqref{eq_mixed_sum}, we have 
 \begin{align*}
     \reg (I^t) & = \reg (P^{t}) + \reg (Q)    - 1 \\
     & = \reg (P) + 2 t-3 + (k-1) \reg (P) - (k-2) -1 \\
     & = k \reg (P) -k + 2t - 2.
 \end{align*}
The equality for symbolic powers follows by the same line of proof by applying \cite[Theorem 4.7]{HNTT} instead of \cite[Theorem 1.1]{NV}. The conclusion follows.
\end{proof}

We are now ready to derive the regularity of powers and symbolic powers of general cubic circulant graphs.
\begin{thm}
	Let $G=C_{2n}(a,n)$ with $1 \le a \le n - 1$ be a cubic circulant graph. Let $d = \gcd(a,2n)$. Then we have $\reg (I(G)^{(t)}) = \reg (I(G)^t)$ for all $t \ge 1$ and 
 $$\reg (I(G)^t)  = \begin{cases}
     d \left \lfloor \frac{n}{2d} \right \rfloor + 2t - 1 & \text{ if } \frac{2n}{d} \text{ is even and }  \frac{n}{d} \not \equiv 1 \pmod 4,\\
     d \left \lfloor \frac{n}{2d} \right \rfloor + d + 2t-2 & \text{ if } \frac{2n}{d} \text{ is even and } \frac{n}{d} \equiv 1 \pmod 4, \\
     \frac{d}{2} \left \lfloor \frac{n}{d} \right \rfloor + 2t - 1 & \text{ if } \frac{2n}{d} \text{ is odd and } \frac{2n}{d} \not \equiv 3 \pmod 4,\\
     \frac{d}{2} \left \lfloor \frac{n}{d} \right \rfloor + \frac{d}{2} + 2t-2 & \text{ if } \frac{2n}{d} \text{ is odd, } \frac{2n}{d} \equiv 3 \pmod 4 \text{ and } \frac{2n}{d} \neq 3, \\  
     \frac{d}{2} + 2t - 1 & \text{ if } \frac{2n}{d} = 3,
 \end{cases}$$
 for all $t \ge 2$.
\end{thm}

\begin{proof}
    Assume that $G = kG_1$ is the disjoint union of $k$ copies of connected cubic circulant graph $G_1 = C_{2m}(a,m)$ for some $a \in \{1,2\}$. By Theorem \ref{thm_main}, Proposition \ref{prop1}, Lemma \ref{lem_mixed_sum_1}, and Lemma \ref{lem_mixed_sum_2} we have 
    $$\reg (I(G)^t) = \reg (I(G)^{(t)}) = \begin{cases}
        k \im (G_1) + 2t - 1 & \text{ if } \reg (I(G_1)) = \im(G_1) + 1,\\
        k \im(G_1) + k + 2t - 2 & \text{ if } \reg (I(G_1)) = \im(G_1)+2,
    \end{cases}$$
    for all $t \ge 2$. The conclusion then follows from Lemma \ref{lem_structure_cubic} and Lemma \ref{lem_induced_matching_cubic}.
\end{proof}
\section{Acknowledgment} 
We thank Prof. Nguyen Cong Minh for many helpful comments and suggestions regarding the regularity of symbolic powers of cubic circulant graphs. A part of this work was completed during the first two authors' visit to the Vietnam Institute for Advanced Study in Mathematics (VIASM). They wish to express their gratitude for VIASM's support and hospitality.

  \section*{Declarations}
  
  
  \subsection*{\bf Ethical approval}\ \\[-0.3cm]
  
  \noindent Not applicable.
  
  
  \subsection*{\bf Competing interests}\ \\[-0.3cm]
  
  \noindent The authors declare that there are no conflict of interests.
  
  
  \subsection*{\bf Authors' contributions}\ \\[-0.3cm]
  
  \noindent All authors have contributed equally to this work.
  

  
  \subsection*{\bf Availability of data and materials}\ \\[-0.3cm]
  
  \noindent Data sharing is not applicable to this work as no datasets were generated or analysed during the current study.

\end{document}